\documentclass[10pt,sumlimits,namelimits]{article}
\newcommand{\nc}{\newcommand}
\newcommand{\rnc}{\renewcommand}
\usepackage{amsthm}
\usepackage{amsmath}
\usepackage{amsfonts}
\usepackage{amssymb}
\usepackage{color}
\usepackage[colorlinks=true,urlcolor=blue]{hyperref}
\usepackage[applemac]{inputenc}

\nc{\N}{\mathbb{N}}
\nc{\Z}{\mathbb{Z}}
\nc{\D}{\mathbb{D}}
\nc{\Q}{\mathbb{Q}}
\nc{\R}{\mathbb{R}}
\nc{\C}{\mathbb{C}}
\nc{\vphi}{\varphi}
\nc{\dsp}{\displaystyle}
\nc{\ovl}{\overline}
\nc{\udl}{\underline}
\nc{\vlim}{\lim\limits}
\nc{\vlimsup}{\limsup\limits}
\nc{\vliminf}{\liminf\limits}
\nc{\vsup}{\sup\limits}
\nc{\vinf}{\inf\limits}
\nc{\vint}{\int\limits}
\nc{\tends}{\longrightarrow}
\nc{\loc}{{\rm loc}}
\rnc{\le}{\leqslant}
\rnc{\ge}{\geqslant}
\rnc{\Re}{{\rm Re}}
\rnc{\Im}{{\rm Im}}

\numberwithin{equation}{section}
\rnc{\theequation}{\thesection.\arabic{equation}}

\newtheorem{thm}{Theorem}[section]
\newtheorem{prop}[thm]{Proposition}

\newtheorem{lem}[thm]{Lemma}

\theoremstyle{definition}
\newtheorem{rmk}[thm]{Remark}
\newtheorem{defi}[thm]{Definition}

\newenvironment{proof*}{\noindent{\bf Proof.}}{\qed}
\newenvironment{vproof}[1]{\noindent{\bf Proof #1}}{\qed}

\title{\Huge \sc Convergence to Scattering States in the Nonlinear Schrödinger Equation}
\author{\sc Pascal Bégout}
\date{}

\begin{document}

\maketitle

\begin{center}
Laboratoire d'Analyse Numérique \\
Université Pierre et Marie Curie \\
Boîte Courrier 187 \\
4, place Jussieu 75252 Paris Cedex 05, FRANCE \bigskip \\
{\footnotesize e-mail\:: }\htmladdnormallink{{\footnotesize\udl{\tt{begout@ann.jussieu.fr}}}}
{mailto:begout@ann.jussieu.fr}
\end{center}

\begin{abstract}
In this paper, we consider global solutions of the following nonlinear Schrödinger equation $iu_t+\Delta u+\lambda|u|^\alpha u =
0,$ in $\R^N,$ with $\lambda\in\R,$ $\alpha\in(0,\frac{4}{N-2})$ $(\alpha\in(0,\infty)$ if $N=1)$ and \linebreak $u(0)\in X\equiv
H^1(\R^N)\cap L^2(|x|^2;dx).$ We show that, under suitable conditions, if the solution $u$ satisfies $e^{-it\Delta}u(t)-u_
\pm\to0$ in $X$ as $t\to\pm\infty$ then $u(t)-e^{it\Delta}u_\pm\to0$ in $X$ as $t\to\pm\infty.$ We  also study the converse. 
Finally, we estimate $|\:\|u(t)\|_X-\|e^{it\Delta}u_\pm\|_X\:|$ under some less restrictive assumptions. 
\end{abstract}

\baselineskip .7cm

\section{Introduction and notations}
\label{introduction}

{\let\thefootnote\relax\footnotetext{2000 Mathematics Subject Classification: 35Q55 (35B40, 35B45, 35P25)}}

We consider the following Cauchy problem,
\begin{gather}
 \left\{
  \begin{split}
   \label{nls}
    i\frac{\partial u}{\partial t}+\Delta u+\lambda|u|^\alpha u & = 0,\; (t,x)\in(-T_*,T^*)\times \R^N, \\
                                                           u(0) & = \vphi, \mbox{ in }\R^N,
  \end{split}
 \right.
\end{gather}
where $\lambda\in \R,$ $0\le\alpha<\dfrac{4}{N-2}$ $(0\le\alpha<\infty$ if $N=1)$ and $\vphi$ a given initial data.

It is well-known that if $\lambda<0,$ $\alpha>\dfrac{4}{N}$ and $\vphi\in H^1(\R^N),$ then there exists $u_\pm\in H^1(\R^N)$
such that $\vlim_{t\to\pm\infty}\|T(-t)u(t)-u_{\pm}\|_{H^1}=0$ (Ginibre and Velo \cite{MR87i:35171}, Nakanishi
\cite{MR1726753,MR1829982}). Since $(e^{it\Delta})_{t\in\R}$ is an isometry on $H^1(\R^N),$ we also have $\vlim_{t\to\pm\infty}\|
u(t)-T(t)u_{\pm}\|_{H^1}=0.$ Furthermore, if $\alpha>\frac{-(N-2)+\sqrt{N^2+12N+4}}{2N}$ and if $\vphi\in X\equiv H^1(\R^N)\cap
L^2(|x|^2;dx),$ then there exist $u_\pm\in X$ such that $\vlim_{t\to\pm\infty}\|T(-t)u(t)-u_{\pm}\|_{X}=0$
(Tsutsumi~\cite{MR87g:35221}). The same result  holds without assumption on the $\lambda$'s sign if the initial data is small enough in $X$
and if $\alpha>\dfrac{4}{N+2}$ (Cazenave and Weissler \cite{MR93d:35150}). Note  that to have these limits, we  have to make a
necessary assumption on $\alpha$ (Barab \cite{MR86a:35121}, Strauss \cite{297.35062,MR83b:47074a}, Tsutsumi and Yajima
\cite{MR85k:35216}), $\dfrac{2}{N}<\alpha<\dfrac{4}{N-2}$ $(2<\alpha<\infty$ if $N=1).$

The purpose of this paper is to study the asymptotic behavior of $\| u(t)-T(t)u_{\pm}\|_{X}$ under the assumption 
$\vlim_{t\to\pm\infty}\|T(-t)u(t)-u_{\pm}\|_{X}=0,$ and the converse. In the linear case $(i.e.:$ $\lambda=0$) or if the
initial data is $0$, the answer is trivial since $T(-t)u(t)-u_\pm\equiv u(t) - T(t)u_\pm\equiv0,$ for all $t\in\R.$ Since
$(e^{it\Delta})_{t\in\R}$ is an isometry on $H^1(\R^N),$ the equivalence on $H^1(\R^N)$ is trivial. But 
$(e^{it\Delta})_{t\in\R}$ is not an  isometry on $X$ and so it is natural to wonder whether or not we have 
$\vlim_{t\to\pm\infty}\| u(t)-T(t)u_{\pm}\|_{X}=0$ when $\vlim_{t\to\pm\infty}\|T(-t)u(t)-u_{\pm}\|_{X}=0$ and 
conversely. \\
\indent
This paper is organized as follows. In Section \ref{mainresult}, we give the main results. In Section \ref{estimapriori}, we
establish some {\it a priori} estimates. In Section \ref{proof1}, we prove Theorems \ref{commut1}, \ref{borne+}, \ref{borne-} and
Proposition \ref{autreborne}. In Section \ref{proof2}, we prove Theorem \ref{commut2}. \\
\indent
Before closing this section, we give some notations which will be used throughout this paper and we recall some
properties of the solutions of the nonlinear Schrödinger equation. \\
\indent
$\ovl{z}$ is the conjugate of the complex number $z;$
$\Re$ $z$ and $\Im$ $z$ are respectively the real and imaginary part of the complex number $z;$
$\Delta=\sum\limits_{j=1}^N\frac{\partial^2}{\partial x_j^2};$
for $1\le p\le\infty,$ $p'$ is the conjugate of the real number $p$ defined by
$\frac{1}{p}+\frac{1}{p'}=1$ and
$L^p= L^p(\R^N)= L^p(\R^N;\C)$ with norm $\|\: .\:\|_{L^p};$
$H^1= H^1(\R^N)= H^1(\R^N;\C)$ with norm $\|\: .\:\|_{H^1};$
for all $(f,g)\in L^2\times L^2,$ $(f,g)=\Re\vint_{\R^N}f(x)\ovl{g(x)}dx;$
$X=\left\{\psi\in H^1(\R^N;\C);\; \|\psi\|_X<\infty\right\}$ with norm $\|\psi\|_X^2=\|\psi\|_{H^1(\R^N)}^2+
\vint_{\R^N}|x|^2|\psi(x)|^2dx;$
$(T(t))_{t \in \R}$ is the group of isometries $(e^{it\Delta})_{t\in\R}$ generated by $i\Delta$ on $L^2(\R^N;\C);$
$C$ are auxiliary positive constants and $C(a_1, a_2,\dots,a_n)$ indicates that the constant $C$ depends only
on parameters $\:a_1,a_2,\dots, a_n$ and that the dependence is continuous. \\
\indent
It is well-known that for every $\vphi\in X,$ (\ref{nls}) has a unique solution $u\in C((-T_*,T^*);X)$ which 
satisfies the  conservation of charge and energy, that is for all $t\in(-T_*,T^*),$ $\|u(t)\|_{L^2}=\|\vphi\|_{L^2}$
and $E(u(t))=E(\vphi)\stackrel{\rm def}{=}\frac{1}{2}\|\nabla\vphi\|_{L^2}^2-\frac{\lambda}{\alpha+2}\|\vphi\|_{L^{\alpha+2}}
^{\alpha+2}.$ Moreover, if $\lambda\le0,$ if $\alpha<\frac{4}{N}$ or if $\|\vphi\|_{H^1}$ is small  enough then
$T^*=T_*=\infty$ and $\|u\|_{L^\infty(\R;H^1)}<\infty$ (see for example Cazenave \cite{caz1}, Ginibre and Velo 
\cite{MR82c:35059,MR82c:35057,MR82c:35058,MR87b:35150}, Kato \cite{MR88f:35133}).

\begin{defi}
\label{padef}
We say that $(q,r)$ is an {\it admissible pair} if the following holds. \medskip \\
$
\begin{array}{rl}
 (\rm i) & 2\le r\le\frac{2N}{N-2}\; (2\le r<\infty\; $if$\; N=2,\; 2\le r\le\infty\mbox\; $if$\; N=1), \medskip \\
(\rm ii) & \frac{2}{q} = N\left(\frac{1}{2} - \frac{1}{r}\right).
\end{array}
$
\\
Note that in this case $2\le q\le\infty$ and $q=\dfrac{4r}{N(r-2)}.$
\end{defi}

\begin{defi}
\label{scatteringdef}
We say that a solution $u\in C((-T_*,T^*);X)$ of (\ref{nls}) has a {\it scattering state} $u_+$ at $+\infty$ (respectively $u_-$
at
$-\infty$) if $T^*=\infty$ and if $u_+\in X$ is such that $\vlim_{t\to\infty}\|T(-t)u(t)-u_+\|_X=0$  (respectively if
$T_*=\infty$  and if $u_-\in X$ is such that $\vlim_{t\to-\infty}\|T(-t)u(t)-u_-\|_X=0$).
\end{defi}
We recall the Strichartz' estimates. Let $I\subseteq\R,$ be an interval, let $t_0\in\ovl I,$ let $(q,r)$ and $(\gamma,\rho)$ be
two admissible pairs, let $\vphi\in L^2(\R^N)$ and let $f\in L^{\gamma'}(I;L^{\rho'}(\R^N)).$ Then the following integral
equation defined for all $t\in I,$ $u(t)=T(t)\vphi+i\dsp\int_{t_0}^tT(t-s)f(s)ds,$ satisfies the following inequality
$
\|u\|_{L^q(I,L^r)}\le C_0\|\vphi\|_{L^2}+C_1\|f\|_{L^{\gamma'}(I;L^{\rho'})},
$
where $C_0=C_0(N,r)$ and $C_1=C_1(N,r,\rho).$ For more details, see Keel and Tao \cite{MR1646048}.

\section{The main results}
\label{mainresult}

\begin{thm}
\label{commut1}
Let $\lambda\neq0,$ $\dfrac{2}{N}<\alpha<\dfrac{4}{N-2}$ $(2<\alpha<\infty$ if $N=1),$ $\vphi\in X$ and let $u$ be the solution of
$(\ref{nls})$ such that $u(0)=\vphi.$ We assume that $u$ has a scattering state $u_\pm$ at $\pm\infty$ $($see Definition
$\ref{scatteringdef}).$ Then the following holds.
\begin{enumerate}
 \item
 \label{commut11}
  \begin{enumerate}
   \item
   \label{commut111}
     If $N\le2$ and if $\alpha>\dfrac{4}{N}$ then $\vlim_{t\to\pm\infty}\|u(t)-T(t)u_\pm\|_X=0.$
   \item
   \label{commut112}
     If $3\le N\le5$ and if $\alpha>\dfrac{8}{N+2}$ then $\vlim_{t\to\pm\infty}\|u(t)-T(t)u_\pm\|_X=0.$
  \end{enumerate}
 \item
 \label{commut12}
     If $N=1$ and $\alpha=4$ or if $3\le N\le5$ and $\alpha=\dfrac{8}{N+2}$ then we have,
       $$
       \sup_{t\ge0}\|u(t)-T(t)u_+\|_X<\infty\quad and \quad\sup_{t\le0}\|u(t)-T(t)u_-\|_X<\infty.
       $$
\end{enumerate}
\end{thm}

\begin{rmk}
Remark that in Theorem \ref{commut1}, no hypothesis on the $\lambda'$ s sign is made.
\end{rmk}

\begin{rmk}
$N\in\{3,4,5\}\Longrightarrow\dfrac{4}{N}<\dfrac{8}{N+2}<\dfrac{6}{N}<\dfrac{4}{N-2}.$\\
$N\ge6\Longrightarrow\dfrac{4}{N-2}\le\dfrac{8}{N+2}.$
\end{rmk}

Despite the fact we do not know if $\vlim_{t\to\pm\infty}\|u(t)-T(t)u_\pm\|_X=0$ when $\alpha\le\dfrac{8}{N+2}$
$(\alpha\le4/N$ if $N\le2)$ or when $N\ge6,$ we can give an estimate of the difference of the norms, as shows the following
theorem, without any restriction on the dimension space $N$ and on $\alpha$ (except $\alpha>\frac{2}{N}).$ Since under the
scattering state assumption we always have $\vlim_{t\to\pm\infty}\|u(t)-T(t)u_\pm\|_{H^1}=0,$ it is sufficient to estimate
$|\:\|xu(t)\|_{L^2}-\|xT(t)u_\pm\|_{L^2}|$ as $t\tends\pm\infty.$
\medskip

\begin{thm}
\label{borne+}
Let $\lambda<0,$ $\dfrac{2}{N}<\alpha<\dfrac{4}{N-2}$ $(2<\alpha<\infty$ if $N=1),$ $\vphi\in X$ and let $u$ be the associated
solution of $(\ref{nls}).$ Assume that $u$ has a scattering state $u_\pm$ at $\pm\infty$ $($see
Definition $\ref{scatteringdef}).$ We define for all $t\in\R,$ ${\cal A}_\pm(t)=\|xu(t)\|_{L^2}-\|xT(t)u_\pm\|_{L^2}$ and
$h(t)=\|xu(t)\|^2_{L^2}.$ Then we have the following result.
$$
\vsup_{t\ge 0}|\:\|u(t)\|_X - \| T(t)u_+\|_X\:| <\infty\quad and \quad\vsup_{t\le 0}|\:\|u(t)\|_X - \| T(t)u_-\|_X\:| <\infty,
$$
with the following estimates.
\begin{enumerate}
\item
\label{borne+1}
If $\alpha <\dfrac{4}{N}$ then $-\dfrac{C}{\|\nabla u_\pm\|_{L^2}} \le \vliminf_{t\to\pm\infty}{\cal A}_\pm(t) \le
\vlimsup_{t\to\pm\infty}{\cal A}_\pm(t) \le \pm\dfrac{h'(0)+4(xu_\pm,i\nabla u_\pm)}{4\|\nabla u_\pm\|_{L^2}}.$
\item
\label{borne+2}
If $\alpha >\dfrac{4}{N}$ then $\pm\dfrac{h'(0)+4(xu_\pm,i\nabla u_\pm)}{4\|\nabla u_\pm\|_{L^2}} \le
\vliminf_{t\to\pm\infty}{\cal A}_\pm(t) \le \vlimsup_{t\to\pm\infty}{\cal A}_\pm(t) \le \dfrac{C}{\|\nabla u_\pm\|_{L^2}}.$
\item
\label{borne+3}
If $\alpha =\dfrac{4}{N}$ then $\vlim_{t\to\pm\infty}{\cal A}_\pm(t)=\pm\dfrac{h'(0)+4(xu_\pm,i\nabla u_\pm)}
{4\|\nabla u_\pm\|_{L^2}}.$
\end{enumerate}
Furthermore, $h'(0)=4\Im\dsp{\vint_{\R^N}}\ovl{\vphi(x)}x.\nabla\vphi(x)dx$ and
$C=C(\vsup_{t\in\R}\|T(-t)u(t)\|_X,N,\alpha,\lambda).$
\end{thm}
\medskip 

\begin{thm}
\label{borne-}
Let $\lambda>0,$ $\dfrac{2}{N}<\alpha<\dfrac{4}{N-2}$ $(2<\alpha<\infty$ if $N=1),$ $\vphi\in X$ and let $u$ be the associated
solution of $(\ref{nls}).$ Assume that $u$ has a scattering state $u_\pm$ at $\pm\infty$ $($see Definition
$\ref{scatteringdef}).$ We define for all $t\in\R,$ ${\cal A}_\pm(t)=\|xu(t)\|_{L^2}-\|xT(t)u_\pm\|_{L^2}$ and
$h(t)=\|xu(t)\|^2_{L^2}.$ Then we have the following result.
$$
\vsup_{t\ge 0}|\:\|u(t)\|_X - \| T(t)u_+\|_X\:| <\infty\quad and \quad\vsup_{t\le 0}|\:\|u(t)\|_X - \| T(t)u_-\|_X\:| <\infty,
$$
with the following estimates.
\begin{enumerate}
\item
\label{borne-1}
If $\alpha <\dfrac{4}{N}$ then $\pm\dfrac{h'(0)+4(xu_\pm,i\nabla u_\pm)}{4\|\nabla u_\pm\|_{L^2}} \le
\vliminf_{t\to\pm\infty}{\cal A}_\pm(t) \le \vlimsup_{t\to\pm\infty}{\cal A}_\pm(t) \le \dfrac{C}{\|\nabla u_\pm\|_{L^2}}.$
\item
\label{borne-2}
If $\alpha >\dfrac{4}{N}$ then $-\dfrac{C}{\|\nabla u_\pm\|_{L^2}} \le \vliminf_{t\to\pm\infty}{\cal A}_\pm(t) \le
\vlimsup_{t\to\pm\infty}{\cal A}_\pm(t) \le \pm\dfrac{h'(0)+4(xu_\pm,i\nabla u_\pm)}{4\|\nabla u_\pm\|_{L^2}}.$
\item
\label{borne-3}
If $\alpha =\dfrac{4}{N}$ then $\vlim_{t\to\pm\infty}{\cal A}_\pm(t)=\pm\dfrac{h'(0)+4(xu_\pm,i\nabla u_\pm)}
{4\|\nabla u_\pm\|_{L^2}}.$
\end{enumerate}
Furthermore, $h'(0)=4\Im\dsp{\vint_{\R^N}}\ovl{\vphi(x)}x.\nabla\vphi(x)dx$ and
$C=C(\vsup_{t\in\R}\|T(-t)u(t)\|_X,N,\alpha,\lambda).$
\end{thm}
\medskip 

\begin{rmk}
\label{5.43}
By Theorem \ref{commut1} and \ref{borne+2} of Theorems \ref{borne+} and \ref{borne-}, if $\alpha>\dfrac{4}{N}$ when
$N\le2$ or if $\alpha>\dfrac{8}{N+2}$ when $N\in\{3,4,5\}$, we have
$$
\begin{array}{rl}
\Im\dsp{\vint_{\R^N}}\ovl{u_-(x)}x.\nabla u_-(x)dx\le\Im\dsp{\vint_{\R^N}}\ovl{\vphi(x)}x.\nabla\vphi(x)dx
                                         \le\Im\dsp{\vint_{\R^N}}\ovl{u_+(x)}x.\nabla u_+(x)dx, & \!\!\mbox{if } \lambda<0, \\
\Im\dsp{\vint_{\R^N}}\ovl{u_+(x)}x.\nabla u_+(x)dx\le\Im\dsp{\vint_{\R^N}}\ovl{\vphi(x)}x.\nabla\vphi(x)dx
                                         \le\Im\dsp{\vint_{\R^N}}\ovl{u_-(x)}x.\nabla u_-(x)dx, & \!\!\mbox{if } \lambda>0.
\end{array}
$$
\end{rmk}

\begin{rmk}
\label{5.45}
If $\vphi\in H^1(\R^N;\C)$ satisfies $\vphi\equiv a\psi$ with $\psi\in H^1(\R^N;\R)$ and $a\in\C,$ then we have 
$h'(0)\equiv\dfrac{d}{dt}\|xu(t)\|^2_{L^2|t=0} \equiv 4\Im\dsp{\vint_{\R^N}}\ovl{\vphi(x)}x.\nabla\vphi(x)dx=0.$
\end{rmk}

The following proposition offers others estimates.

\begin{prop}
\label{autreborne}
Let $\lambda\neq0,$ $\dfrac{2}{N}<\alpha<\dfrac{4}{N-2}$ $(2<\alpha<\infty$ if $N=1),$ $\vphi\in X$ and let $u$ be the solution
of $(\ref{nls})$ such that $u(0)=\vphi.$ Assume that $u$ has a scattering state $u_\pm$ at $\pm\infty.$ Then the following
estimates hold.
\begin{enumerate}
\item
\label{autreborne1}
If $\lambda<0,$ $\vlimsup_{t\to\pm\infty}(\|xu(t)\|_{L^2}-\|xT(t)u_\pm\|_{L^2})\le
\dfrac{\|xu_\pm\|_{L^2}\|\nabla u_\pm\|_{L^2}\pm(xu_\pm,i\nabla u_\pm)}{\|\nabla u_\pm\|_{L^2}}.$
\item
\label{autreborne2}
If $\lambda>0,$ $\vliminf_{t\to\pm\infty}(\|xu(t)\|_{L^2}-\|xT(t)u_\pm\|_{L^2})\ge
-\dfrac{\|xu_\pm\|_{L^2}\|\nabla u_\pm\|_{L^2}\mp(xu_\pm,i\nabla u_+)}{\|\nabla u_\pm\|_{L^2}}.$
\end{enumerate}
\end{prop}

\begin{rmk}
\label{5.47}
By \ref{borne+3} of Theorems \ref{borne+} and \ref{borne-} and by Proposition \ref{autreborne}, if $\alpha=\dfrac{4}{N}$ then we
have,
$$
\begin{array}{rl}
-\|xu_-\|_{L^2}\|\nabla u_-\|_{L^2}\le\dfrac{1}{4}\dfrac{d}{dt}\|xu(t)\|^2_{L^2|t=0}\le\|xu_+\|_{L^2}\|\nabla u_+\|_{L^2}, &
                                                                                       \!\!\mbox{if } \lambda<0, \medskip \\
-\|xu_+\|_{L^2}\|\nabla u_+\|_{L^2}\le\dfrac{1}{4}\dfrac{d}{dt}\|xu(t)\|^2_{L^2|t=0}\le\|xu_-\|_{L^2}\|\nabla u_-\|_{L^2}, &
                                                                                       \!\!\mbox{if } \lambda>0.
\end{array}
$$
\end{rmk}
\medskip

Now we give the result concerning the converse.

\begin{thm}
\label{commut2}
Let $\lambda\neq0,$ $\dfrac{2}{N}<\alpha<\dfrac{4}{N-2}$ $(2<\alpha<\infty$ if $N=1),$ $\vphi\in X$ and $u$ be the associated
solution of $(\ref{nls}).$ Assume that $u$ is global in time and there exists $u_+\in X$ and $u_-\in X$ such that
$\vlim_{t\to\pm\infty}\|u(t)-T(t)u_\pm\|_X=0.$ Let $\alpha_0=\dfrac{-(N-2)+\sqrt{N^2+12N+4}}{2N}.$ Then, we have the
following result.
\begin{enumerate}
\item
\label{commut21}
If $\lambda<0$ and if $\alpha\ge\alpha_0$ $(\alpha>\alpha_0$ if $N=2)$ then $\vlim_{t\to\pm\infty}\|T(-t)u(t)-u_\pm\|_X=0.$
\item
\label{commut22}
If $\lambda>0$ and if $\alpha>\dfrac{4}{N}$ then $\vlim_{t\to\pm\infty}\|T(-t)u(t)-u_\pm\|_X=0.$
\item
\label{commut23}
If $\alpha>\dfrac{4}{N+2}$ and if $\|\vphi\|_X$ is small enough then $\vlim_{t\to\pm\infty}\|T(-t)u(t)-u_\pm\|_X=0.$
\end{enumerate}
\end{thm}
\medskip

\begin{rmk}
\label{rmkcommut2}
Note that in the case \ref{commut22}, no hypothesis on the $\|\vphi\|_X$' s size is made.
\end{rmk}

\begin{rmk}
\label{rmkcommut3}
Assume there exists $u_\pm,v_\pm\in X$ such that $\vlim_{t\to\pm\infty}\|T(-t)u(t)-u_\pm\|_X=0$ and
$\vlim_{t\to\pm\infty}\|u(t)-T(t)v_\pm\|_X=0.$ Then we have, $u_+=v_+$ and $u_-=v_-.$ Indeed, since $X\hookrightarrow L^2(\R^N)$
and $T(t)$ is an isometry on $L^2(\R^N),$ we have $\vlim_{t\to\pm\infty}\|T(-t)u(t)-u_\pm\|_{L^2}
=\vlim_{t\to\pm\infty}\|T(-t)u(t)-v_\pm\|_{L^2}=0.$ Hence the result.
\end{rmk}

\begin{rmk}
$\alpha_0\in\left(\frac{4}{N+2},\frac{4}{N}\right)$ $(\alpha_0\in\left(\frac{2}{N},\frac{4}{N}\right)$ if $N=1).$
\end{rmk}

\section{A priori estimates}
\label{estimapriori}

Throughout this section, we make the following assumptions.
\begin{gather}
 \begin{cases}
  \label{hypo}
   \lambda\neq0,\; \dfrac{2}{N}<\alpha<\dfrac{4}{N-2}\; (2<\alpha<\infty) \mbox{ if } N=1),\; \vphi\in X,\; u\in C(\R;X)
                   \mbox{ is the} \\
   \mbox{associated solution of } (\ref{nls}) \mbox{ and has a scattering state } u_\pm\in X \mbox{ at } \pm\infty.
 \end{cases}
\end{gather}
We define the following real.
\begin{gather}
\label{gamma*}
\gamma^*=
 \begin{cases}
  \dfrac{\alpha-2}{2},      & \text{if }N=1, \medskip \\
  \dfrac{\alpha(N+2)-4}{4}, & \text{if }N\ge2.
 \end{cases}
\end {gather}

\begin{prop}
\label{propdecayX}
Assume $u$ satisfies $(\ref{hypo})$ $($we can suppose instead of $u$ has a scattering state that we only have
$\vsup_{t\in\R}\|T(-t)u(t)\|_X<\infty).$ Let $(q,r)$ be an admissible pair $($see Definition $\ref{padef}).$ Then the following
holds.
\begin{enumerate}
\item
\label{propdecayX1}
For all $t\neq0,$ $\|u(t)\|_{L^r}\le C|t|^{-\frac{2}{q}},$ where $C=C(\vsup_{t\in\R}\|T(-t)u(t)\|_X,N,r).$
\item
\label{propdecayX2}
If furthermore $\alpha>\dfrac{4}{N+2}$ then $u\in L^q(\R;W^{1,r}(\R^N)).$
\end{enumerate}
\end{prop}

\begin{proof*}
We follow the method of  Cazenave \cite{caz1} (see Theorem 7.2.1 and Corollary 7.2.4). We set
$w(t,x)=e^{-i\frac{|x|^2}{4t}}u(t,x)$ and $f(u)=\lambda|u|^\alpha u.$ We already know that for every admissible pair $(q,r),$
$u\in L^q_\loc(\R;W^{1,r}(\R^N))$ (see for example Cazenave \cite{caz1}; Theorem 5.3.1 and Remark 5.3). We only prove the case
$t>0,$ the case $t<0$ following by applying the result for $t>0$ to $\ovl{u(-t)}$ solution of (\ref{nls}) with initial value
$\ovl{\vphi}.$ We proceed in 2 steps. \\
{\bf Step 1.} $\|u(t)\|_{L^r}\le
C(\vsup_{t\in\R}\|T(-t)u(t)\|_X,N,r)|t|^{-\frac{2}{q}},$ for every admissible pair $(q,r)$ and for all $t\neq0.$ \\
We have $\|xT(-t)u(t)\|_{L^2}=\|(x+2it\nabla)u(t)\|_{L^2}\le C.$ Furthermore,
$
(x+2it\nabla)u(t,x)=2ite^{-i\frac{|x|^2}{4t}}\nabla w(t,x).
$
Using the Gagliardo-Nirenberg's inequality, we obtain
\begin{align*}
\|u(t)\|_{L^r}\equiv\|w(t)\|_{L^r} & \le C\|\nabla w(t)\|_{L^2}^{N\left(\frac{1}{2}-\frac{1}{r}\right)}\|
                                           w(t)\|_{L^2}^{1-N\left(\frac{1}{2}-\frac{1}{r}\right)} \\
                                   & \le C \left(\|(x+2it\nabla)u(t)\|_{L^2}|t|^{-1}\right)^{N\left(\frac{1}{2}
                                                                             -\frac{1}{r}\right)} \\
                                   & \le C |t|^{-N\left(\frac{1}{2}-\frac{1}{r}\right)}.
\end{align*}
Hence the result. \\
{\bf Step 2.} $u\in L^q(\R;W^{1,r}(\R^N))$ for every admissible pair $(q,r).$ \\
By the Strichartz' estimates and by Hölder's inequality (applying twice), we have
\begin{gather}
\label{dempropdecayX2}
\|f(u)\|_{L^{q'}((0,\infty);W^{1,r'})}\le C \left( \|u\|_{L^{\frac{q\alpha}{q-2}}((0,1);L^{\frac{r\alpha}{r-2}})}+
\|u\|_{L^{\frac{q\alpha}{q-2}}((1,\infty);L^{\frac{r\alpha}{r-2}})}\right)^\alpha\|u\|_{L^q((0,\infty);W^{1,r})}, \\
\label{dempropdecayX1}
\|u\|_{L^q((S,\infty);W^{1,r})} \le C + C \|u\|_{L^{\frac{q\alpha}{q-2}}((S,\infty);L^{\frac{r\alpha}{r-2}})}
                                                                       ^\alpha\|u\|_{L^q((S,\infty);W^{1,r})},
\end{gather}
for all $S>0$ and for every admissible pair $(q,r).$ \\
{\it Case N$\:\ge$3.} We set $r=\frac{4N}{2N-\alpha(N-2)}.$ Since $\alpha\in\left(0,\frac{4}{N-2}\right),$ we have 
$r\in\left(2,\frac{2N}{N-2}\right).$ So we can take $q$ such that $(q,r)$ is an admissible pair. For this choice of $r,$ we have 
$\frac{r\alpha}{r-2}=\frac{2N}{N-2}$ and $\frac{q}{q-2}=\frac{4}{4-\alpha(N-2)}.$ By (\ref{dempropdecayX1}) and the first step
we  have for all $S>0,$
$$
\|u\|_{L^q((S,\infty);W^{1,r})} \le C + C\left(\dsp\int_S^\infty
t^{-\frac{4\alpha}{4-\alpha(N-2)}}dt\right)^{\frac{q-2}{q}}\|u\|_{L^q((S,\infty);W^{1,r})}.
$$
And $\dfrac{4\alpha}{4-\alpha(N-2)}>1\iff\alpha>\dfrac{4}{N+2}.$ Thus, there exists $S_0>0$ large enough such that 
$$
C\left(\dsp\int_{S_0}^\infty t^{-\frac{4\alpha}{4-\alpha(N-2)}}dt\right)^{\frac{q}{q-2}}\le \frac{1}{2},
$$
and then,
$$
\|u\|_{L^q((S_0,\infty);W^{1,r})} \le 2C.
$$
For this choice of $(q,r),$ we deduce from (\ref{dempropdecayX2}) that $\|f(u)\|_{L^{q'}((0,\infty);W^{1,r'})}<\infty.$ Hence the
result for every admissible pair by the Strichartz' estimates. \\
{\it Case N=2.} Since $\alpha>1$ is fixed, we take $r>2$ sufficiently close to 2 to have $\alpha>\dfrac{2(r-1)}{r}.$ 
So, in particular, $\dfrac{r\alpha}{r-2}>2.$ Moreover, $\dfrac{q}{q-2}=\dfrac{r}{2}$ where $q$ is such that $(q,r)$ is an
admissible  pair. So by  Hölder's inequality (twice), (\ref{dempropdecayX1}) and the first step, we have for all $S>0,$
$$
\|u\|_{L^q((S,\infty);W^{1,r})} \le C+C\left(\dsp\int_S^\infty t^{-\frac{r\alpha-2(r-2)}{2}}dt\right)^{\frac{2}{r}}
                                                                       \|u\|_{L^q((S,\infty);W^{1,r})}.
$$
And $\dfrac{r\alpha-2(r-2)}{2}>1\iff\alpha>\dfrac{2(r-1)}{r}.$ And we conclude exactly as the case $N\ge3.$ \\
{\it Case N=1.} We take (\ref{dempropdecayX2}) with the admissible pair $(\infty,2)$ and apply the first step. So,
\begin{align*}
\|f(u)\|_{L^1((0,\infty);H^1)} & \le C \left( \|u\|_{L^\alpha((0,1);L^\infty)}+\|u\|_{L^\alpha((1,\infty);L^\infty)}
                                      \right)^\alpha\|u\|_{L^\infty((0,\infty);H^1)} \bigskip \\
                               & \le C+C \|u\|_{L^\alpha((1,\infty),L^\infty)}^\alpha
                                      \le C+C \dsp\vint_1^\infty t^{-\frac{\alpha}{2}}dt <\infty.
\end{align*}
Hence the result for every admissible pair by the Strichartz' estimates.
\medskip
\end{proof*}

\begin{rmk}
\label{rmkdecayX}
We set $v(t,x)=(x+2it\nabla)u(t,x).$ In the same way as the above proof, we can show under the assumptions of Proposition
\ref{propdecayX} that if $\alpha>\dfrac{4}{N+2},$ then for every admissible pair $(q,r),$ $v\in L^q(\R;L^r(\R^N))$ (follow the
step 2 of the proof of Corollary 7.2.4 of Cazenave \cite{caz1}; consider separately the three cases $N=1,$ $N=2,$ $N\ge3$ and
replace the admissible pairs therein by those of the proof of Proposition \ref{propdecayX}).
\end{rmk}

\begin{prop}
\label{propvitesseX}
Let $\gamma^*$ be defined by $(\ref{gamma*}).$ Assume $u$ satisfies $(\ref{hypo})$ and $\alpha>\dfrac{4}{N+2}.$ Then, the
following estimates hold.
\begin{enumerate}
\item
\label{propvitesseX1}
If $N=1$ then for all $t\neq0$ we have $\|T(-t)u(t)-u_\pm\|_{H^1}\le C|t|^{-\gamma^*}.$
\item
\label{propvitesseX2}
If $N=2$ then for all $t\neq0$ and for any $\gamma<\gamma^*,$ we have $\|T(-t)u(t)-u_\pm\|_{H^1}\le C|t|^{-\gamma}.$
\item
\label{propvitesseX3}
If $N\ge3$ then for all $t\neq0$ we have $\|T(-t)u(t)-u_\pm\|_{H^1}\le C|t|^{-\gamma^*}.$
\end{enumerate}
\end{prop}
\medskip

\begin{proof*}
Denote $f(u)=\lambda|u|^\alpha u.$ We only prove the case $t>0,$ the case $t<0$ following by applying the result for 
$t>0$ to $v(t)=\ovl{u(-t)}$ solution of (\ref{nls}) with $v(0)=\ovl{\vphi}.$ In this case, $v_+=\ovl{u_-}$
and the result follows. By applying the Strichartz' estimates and Hölder's inequality (twice), we have
$$
\|u(t)-T(t)u_+\|_{H^1}\le C\|f(u)\|_{L^{q'}(t,\infty);W^{1,r'})}
                      \le C\|u\|_{L^{\frac{q\alpha}{q-2}}((t,\infty);L^{\frac{r\alpha}{r-2}})}^\alpha\|u\|_{L^q(\R;W^{1,r})},
$$
for every admissible pair $(q,r)$ and for all $t>0.$ Thus, by Proposition \ref{propdecayX}, \ref{propdecayX2}, we have
\begin{gather*}
\|u(t)-T(t)u_+\|_{H^1}\le C\|u\|_{L^{\frac{q\alpha}{q-2}}((t,\infty);L^{\frac{r\alpha}{r-2}})}^\alpha,
\end{gather*}
for every admissible pair $(q,r)$ and for all $t>0.$ \\
Now, we conclude by the same way than for the step 2 of the proof of Proposition \ref{propdecayX}, using \ref{propdecayX1} of this
proposition, considering separately the three cases $N=1,$ $N=2,$ $N\ge3,$ and using the same admissible pairs. This achieves the
proof.
\end{proof*}

\section{Proof of Theorems \ref{commut1}, \ref{borne+}, \ref{borne-} and Proposition \ref{autreborne}}
\label{proof1}

Throughout this section, we assume that $u$ satisfies (\ref{hypo}). \medskip

\begin{vproof}{of Theorem \ref{commut1}.}
Since $\vlim_{t\to\pm\infty}\| T(-t)u(t)-u_\pm\|_X=0$ and $T(t)$ is an isometry on $H^1$, we have 
$\vlim_{t\to\pm\infty}\|u(t)-T(t)u_\pm\|_{H^1}=0$. Thus, it is sufficient to prove that 
$\vlim_{t\to\pm\infty}\|xu(t)-xT(t)u_\pm\|_{L^2}=0$ to obtain \ref{commut11} and that
$\vsup_{t\ge0}\|xu(t)-xT(t)u_+\|_{L^2}<\infty$  and $\vsup_{t\le 0}\|xu(t)-xT(t)u_-\|_{L^2}<\infty$ to obtain
\ref{commut12}. Suppose that the result is proved for $t>0.$ Then we apply it to $v(t)=\ovl{u(-t)}$ solution of (\ref{nls})
with initial data $v(0)=\ovl\vphi.$ Then $v_+=\ovl{u_-}$ is the scattering state at $+\infty$ of $\ovl{u(-t)}.$ And
using the identity $T(t)\ovl{\psi}=\ovl{T(-t)\psi}$ which holds for all $t\in\R$ and for every $\psi\in L^2,$ we obtain
the result for $t<0.$ So to conclude, it is sufficient to prove that $\vlim_{t\to\infty}\|xu(t)-xT(t)u_+\|_{L^2}=0$ to 
obtain \ref{commut11} and $\vsup_{t\ge 0}\|xu(t)-xT(t)u_+\|_{L^2}<\infty$ to obtain \ref{commut12}. We have
\begin{align*}
xu(t)-xT(t)u_+ & = xu(t)-T(t)xu_++2itT(t)\nabla u_+ \\
               & = xu(t)+2it\nabla u(t)-T(t)xu_++2itT(t)\nabla u_+-2it\nabla u(t) \\
			            & = T(t)\left[(xT(-t)u(t)-xu_+)+2it(\nabla u_+-T(-t)\nabla u(t))\right],
\end{align*}
for all $t>0.$ With Proposition \ref{propvitesseX}, we obtain
\begin{align*}
\|xu(t)-xT(t)u_+\|_{L^2} & \le \|xT(-t)u(t)-xu_+\|_{L^2} + 2t\|T(-t)\nabla u(t)-\nabla u_+\|_{L^2} \\
                         & \le \|xT(-t)u(t)-xu_+\|_{L^2} + Ct^{-(\gamma-1)},
\end{align*}
for all $t>0$ and for all $\gamma\in(0,\gamma^*]$ if $N\not=2$ and for all $\gamma\in(0,\gamma^*)$ if $N=2,$ where $\gamma^*$ is
defined by  (\ref{gamma*}). And by assumption, $\vlim_{t\to\infty}\|xT(-t)u(t)-xu_+\|_{L^2}=0$ and $\gamma^* -1>0$ if and 
only if $\alpha>\frac{8}{N+2}$ if $N\ge3$ and $\alpha>\frac{4}{N}$ if $N\le2.$ Thus, in the above expression, it is sufficient 
to choose $\gamma=\gamma^*$ if $N\not=2,$ and $\gamma\in(0,\gamma^*)$ close enough to $\gamma^*$ if $N=2.$ Hence the result.
\medskip
\end{vproof}

\begin{rmk}
Since we do not have the estimate \ref{propvitesseX2} of Proposition \ref{propvitesseX} for $\gamma=\gamma^*,$ we do not know
whether or not $\vsup_{t\ge0}\|u(t)-T(t)u_+\|_X<\infty$ and $\vsup_{t\le0}\|u(t)-T(t)u_-\|_X<\infty$ when $N=2$ and
$\alpha=\dfrac{4}{N}.$
\end{rmk}

We define the following function $h$ by
\begin{gather}
\label{h}
\forall t\in(-T_*,T^*),\; h(t)=\|xu(t)\|^2_{L^2}.
\end{gather}

\begin{lem}
\label{5.5}
Let $u$ satisfying $(\ref{hypo})$ and let $h$ be defined by $(\ref{h}).$ Then $h\in C^2(\R),$ and we have
$$
\begin{array}{rl}
  \forall t\in\R, &
    \left\{
      \begin{array}{rcl}
        h'(t)  & = & 4\Im\dsp{\vint_{\R^N}}\ovl{u(t,x)}x.\nabla u(t,x)dx,\\
        h''(t) & = & 2N\alpha\|\nabla u_+\|_{L^2}^{2} -2(N\alpha-4)\|\nabla u(t) \|_{L^2}^{2}.
      \end{array}
    \right.
\end{array}
$$
Furthermore, if $\lambda<0,$
\begin{gather}
\label{5.51}
\forall t\in\R,\; \|\nabla u(t)\|_{L^2}\le\| \nabla u_+ \|_{L^2},
\end{gather}
 and if $\lambda>0,$
\begin{gather}
\label{5.51'}
\forall t\in\R,\; \|\nabla u(t)\|_{L^2}\ge\| \nabla u_+ \|_{L^2}.
\end{gather}
\end{lem}

\begin{proof*}
See Ginibre and Velo \cite{MR82c:35058} or Cazenave \cite{caz1}, Proposition 6.4.2 to have $h\in C^2(\R)$, the expression of 
$h'$ and $\forall t\in\R$, $h''(t)=4N\alpha E(\vphi)-2(N\alpha-4)\|\nabla u(t)\|_{L^2}^{2}$. Furthermore, using the conservation
of energy and $\vlim_{t\to\pm\infty}\|u(t)\|_{L^{\alpha +2}}=0$ (Proposition \ref{propdecayX}), we obtain $\|\nabla u_+\|_{L^2}^2
=2E(\vphi).$ Which gives, with the above identity, the desired expression of $h''.$ Finally, with the equality $\|\nabla
u_+\|_{L^2}^2 =2E(\vphi)$ and the conservation of energy, we obtain (\ref{5.51}) if $\lambda<0$ and (\ref{5.51'}) if $\lambda>0.$
\medskip
\end{proof*}

Now, we establish 2 lemmas which will be used to prove Theorem \ref{borne+}.

\begin{lem}
\label{5.7}
Let $u$ satisfying $(\ref{hypo})$ and let $h$ be defined by $(\ref{h}).$ Assume that $\lambda<0.$ Then the following
holds. \medskip \\
$
\begin{array}{rl}
\text{\rm (i)}   & If\;\; \alpha<\dfrac{4}{N}\;\; then\; \vlimsup_{t\to\infty}(\|xu(t)\|_{L^2}-\|xT(t)u_+\|_{L^2})
                                 \le\dfrac{h'(0)+4(xu_+,i\nabla u_+)}{4\|\nabla u_+\|_{L^2}}. \medskip \\
\text{\rm (ii)}  & If\;\; \alpha>\dfrac{4}{N}\;\; then\; \vliminf_{t\to\infty}(\|xu(t)\|_{L^2}-\|xT(t)u_+\|_{L^2})
                                 \ge\dfrac{h'(0)+4(xu_+,i\nabla u_+)}{4\|\nabla u_+\|_{L^2}}. \medskip \\
\text{\rm (iii)} & If\;\; \alpha=\dfrac{4}{N}\;\; then\;\vlim_{t\to\infty}(\|xu(t)\|_{L^2}-\|xT(t)u_+\|_{L^2})
                                 =\dfrac{h'(0)+4(xu_+,i\nabla u_+)}{4\| \nabla u_+\|_{L^2}}.
\end{array}
$
\end{lem}
\medskip

\begin{proof*} We proceed in  4 steps.\\
{\bf Step 1.}
$($a$)\quad$If $\alpha <\dfrac{4}{N}$ then $\forall t\ge 0$, $h'(t)\ge h'(0)+2N\alpha\| \nabla u_+ \|_{L^2}^{2}t.$\\
$($b$)\quad$If $\alpha >\dfrac{4}{N}$ then $\forall t\ge 0$, $h'(t)\ge h'(0)+8\| \nabla u_+ \|_{L^2}^{2}t.$\\
$($c$)\quad$If $\alpha =\dfrac{4}{N}$ then $\forall t\in\R$, $h'(t)=h'(0)+8\| \nabla u_+ \|_{L^2}^{2}t.$
\\
We integrate between 0 and $t\ge 0$ the function $h''$ of Lemma \ref{5.5}.\\
$\alpha\le\dfrac{4}{N}\Longrightarrow -2(N\alpha -4)\ge0\Longrightarrow$ (a).\\
$\alpha\ge\dfrac{4}{N}\Longrightarrow -2(N\alpha -4)\le0$ with $(\ref{5.51})\Longrightarrow$ (b).\\
$\alpha=\dfrac{4}{N}\Longrightarrow -2(N\alpha -4)=0\Longrightarrow$(c).
\\
{\bf Step 2.}
$($a$)\quad$If $\alpha <\dfrac{4}{N}$ then $\forall t\ge 0$, $h'(t)\le h'(0)+8\| \nabla u_+ \|_{L^2}^{2}t.$\\
$($b$)\quad$If $\alpha >\dfrac{4}{N}$ then $\forall t\ge 0$, $h'(t)\le h'(0)+2N\alpha\| \nabla u_+ \|_{L^2}^{2}t.$
\\
We integrate between 0 and $t\ge 0$ the function $h''$ of Lemma \ref{5.5}.\\
$\alpha <\dfrac{4}{N}\Longrightarrow -2(N\alpha -4)>0$ and $(\ref{5.51}) \Longrightarrow$(a).\\
$\alpha >\dfrac{4}{N}\Longrightarrow$(b).
\\
{\bf Step 3.}
$($a$)\quad$If $\alpha <\dfrac{4}{N}$ then $\forall t\ge 0,$
$$
\|x\vphi\|_{L^2}^2+h'(0)t+N\alpha\|\nabla u_+\|_{L^2}^2t^2\le h(t) \le \|x\vphi\|_{L^2}^{2}+h'(0)t+4\|\nabla u_+\|_{L^2}^2t^2.
$$
$($b$)\quad$If $\alpha >\dfrac{4}{N}$ then $\forall t\ge 0,$
$$
\|x\vphi\|_{L^2}^2+h'(0)t+4\|\nabla u_+\|_{L^2}^2t^2\le h(t)\le\|x\vphi\|_{L^2}^2+h'(0)t+N\alpha\|\nabla u_+\|_{L^2}^2t^2.
$$
$($c$)\quad$If $\alpha=\dfrac{4}{N}$ then $\forall t\in\R,$ $h(t)=\|x\vphi\|_{L^2}^2+h'(0)t+4\|\nabla u_+\|_{L^2}^2t^2.$
\\
It is sufficient to integrate between 0 and $t\ge0$ the formulas of steps 1 and 2 to obtain the step 3.
\\
{\bf Step 4.} Conclusion.\\
We set\:: $g(t)=\sqrt{\|x\vphi\|_{L^2}^2 + h'(0)t + 4 \| \nabla u_+\|_{L^2}^2t^2},$ $t>0$ large enough. \\
Then for $t>0$ large enough, we have the following asymptotic development:
\begin{gather}
\label{5.02}
g(t)=2\|\nabla u_+\|_{L^2}t+\dfrac{h'(0)}{4\|\nabla u_+\|_{L^2}}+\dfrac{\|x\vphi\|_{L^2}^2}{4\|\nabla
u_+\|_{L^2}t}+o\left(\dfrac{1}{t}\right).
\end{gather}
In the same way, for $t>0$ large enough, we have :
\begin{gather}
\label{5.03}
\|xT(t)u_+\|_{L^2}=2\|\nabla u_+\|_{L^2}t-\dfrac{(xu_+,i\nabla u_+)}{\|\nabla u_+\|_{L^2}}+\dfrac{\|xu_+\|_{L^2}^2}{4\|\nabla
                    u_+\|_{L^2}t}+o\left(\dfrac{1}{t}\right).
\end{gather}
And, applying the step 3 (a), (\ref{5.02}) and (\ref{5.03}) and taking $\vlimsup_{t\to\infty}$, we get (i).
Indeed,
$$
\begin{array}{rl}
    & \|xu(t)\|_{L^2}-\|xT(t)u_+\|_{L^2} \medskip \\
\le & 2\|\nabla u_+\|_{L^2}t+\dfrac{h'(0)}{4\|\nabla u_+\|_{L^2}}+\dfrac{\|x\vphi\|_{L^2}^2}{4\|\nabla u_+\|_{L^2}t}
                      -2\|\nabla u_+\|_{L^2}t+\dfrac{(xu_+,i\nabla u_+)}{\|\nabla u_+\|_{L^2}} - \dfrac{\| xu_+\|_{L^2}^2}
                      {4\|\nabla u_+\|_{L^2} t}+o\left(\dfrac{1}{t}\right) \medskip \\
  = & \dfrac{h'(0)}{4\|\nabla u_+\|_{L^2}}+\dfrac{\|x\vphi\|_{L^2}^2}{4\|\nabla u_+\|_{L^2}t}
                      +\dfrac{(xu_+,i\nabla u_+)}{\|\nabla u_+\|_{L^2}} - \dfrac{\| xu_+\|_{L^2}^2}
                      {4\|\nabla u_+\|_{L^2} t}+o\left(\dfrac{1}{t}\right).
\end{array}
$$
Hence (i) by taking $\vlimsup_{t\to\infty}$ in the above expression. \\
By applying the step 3 (b), (\ref{5.02}) and (\ref{5.03}) and taking $\vliminf_{t\to\infty}$, we get (ii) by the same way.
\\
By applying the step 3 (c), (\ref{5.02}) and (\ref{5.03}) and letting $t\tends\infty,$ we get (iii) by the same way.
\end{proof*}

\begin{lem}
\label{5.8}
Let $u$ satisfying $(\ref{hypo}).$ Assume that $\lambda<0.$ Then the following
holds. \medskip \\
$
\begin{array}{rl}
 ({\rm i}) & If\;\; \alpha <\dfrac{4}{N}\;\; then\; \vliminf_{t\to\infty}(\|xu(t)\|_{L^2}-\|xT(t)u_+\|_{L^2})
             \ge -\dfrac{C}{\|\nabla u_\pm\|_{L^2}}. \medskip \\
({\rm ii}) & If\;\; \alpha >\dfrac{4}{N}\;\; then\; \vlimsup_{t\to\infty}(\|xu(t)\|_{L^2}-\|xT(t)u_+\|_{L^2})
             \le \dfrac{C}{\|\nabla u_\pm\|_{L^2}}.
\end{array}
$
\medskip \\
Furthermore, $C=C(\vsup_{t\in\R}\|T(-t)u(t)\|_X,N,\alpha,\lambda).$
\end{lem}
\medskip

\begin{proof*} We proceed in 2 steps. Let $h$ be defined by (\ref{h}). \\
{\bf Step 1.}
$\forall t\ge 1$, $-\dsp\vint_{1}^{t}\dsp\vint_{1}^{s}\|\nabla u(\sigma)\|_{L^2}^{2}d\sigma ds
\le C+Ct^{-\frac{N\alpha -4}{2}}+Ct+\|\nabla u_+\|^{2}t-\dfrac{1}{2}\|\nabla u_+\|_{L^2}^2t^2,$ where
$C=C(\vsup_{t\in\R}\|T(-t)u(t)\|_X,N,\alpha,\lambda).$ \\
By Proposition \ref{propdecayX}, we have $\|u(\sigma)\|_{L^{\alpha +2}}^{\alpha+2}\le
C(\vsup_{t\in\R}\|T(-t)u(t)\|_X,N,\alpha)\sigma^{-\frac{N\alpha}{2}},$ for all $\sigma >0.$ With the  conservation of energy
and the formula $\| \nabla u_+\|_{L^2}^2=2E(\vphi)$, we deduce that for all $\sigma >0,$
$$
\| \nabla u_+\|_{L^2}^2 - \| \nabla u(\sigma)\|_{L^2}^2\le
C(\vsup_{t\in\R}\|T(-t)u(t)\|_X,N,\alpha,\lambda)\sigma^{-\frac{N\alpha}{2}}
$$
(since $\lambda<0).$ Integrating this expression over $[1,t]\times [1,s]$, we obtain the desired result.
\\
{\bf Step 2.} Conclusion. \\
(i)\quad Lemma \ref{5.5} implies that for every $t\ge0,$
\begin{gather}
\label{5.81}
h(t)=\| x\vphi\|_{L^2}^2+h'(0)t + N\alpha\|\nabla u_+\|_{L^2}^2t^2-2(N\alpha-4)\vint_0^t\vint_0^s\| 
\nabla u(\sigma)\|_{L^2}^2 d\sigma ds.
\end{gather}
By (\ref{5.81}) and step 1, we obtain,
$$\forall t\ge 1,\;
h(t)\ge C+h'(0)t+N\alpha\|\nabla u_+\|_{L^2}^2t^2+Ct^{\frac{4-N\alpha}{2}}-Ct+2(N\alpha -4)\|\nabla u_+\|_{L^2}^2t-(N\alpha -
4)\|\nabla u_+\|^2t^2.
$$
And so for all $t>1,$ $h(t)\ge C+Ct^{\frac{4-N\alpha}{2}}+(h'(0)-2(4-N\alpha)\|\nabla u_+\|_{L^2}^2-C)t+4\|\nabla
u_+\|_{L^2}^2t^2.$
By an asymptotic development on this last inequality, we obtain
$$
\|xu(t)\|_{L^2}\ge2\|\nabla u_+\|_{L^2}t+Ct^{-1}+Ct^{-\frac{N\alpha-2}{2}} \\
+\frac{h'(0)-2(4-N\alpha)\|\nabla u_+\|_{L^2}^2-C}{4\|\nabla u_+\|_{L^2}}+o\left(\frac{1}{t}\right),
$$
for all $t>0$ large enough. From this last expression and (\ref{5.03}), we obtain (i) (see the end of the proof of (i) of Lemma
\ref{5.7}).\\ (ii)\quad From (\ref{5.81}) and step 1, we obtain, for all $t>1,$
$$
h(t)\le C+f'(0)t+N\alpha\|\nabla u_+\|_{L^2}^{2}t^2+Ct^{\frac{4-N\alpha}{2}}
+Ct+2(N\alpha -4)\|\nabla u_+\|_{L^2}^2t - (N\alpha-4)\|\nabla u_+\|_{L^2}^2t^2.
$$
And so, for all $t>1,$ $h(t)\le C+Ct^{\frac{4-N\alpha}{2}}+(h'(0)+2(N\alpha-4)\|\nabla u_+\|_{L^2}^2+C)t+4\|\nabla
u_+\|_{L^2}^2t^2.$
By an asymptotic development on this last inequality, we obtain
$$
\|xu(t)\|_{L^2}\le2\|\nabla u_+\|_{L^2}t+Ct^{-1}+Ct^{-\frac{N\alpha-2}{2}}
+\frac{h'(0)-2(4-N\alpha)\|\nabla u_+\|_{L^2}^2+C}{4\|\nabla u_+\|_{L^2}}+o\left(\frac{1}{t}\right),
$$
for all $t>0$ large enough. With this last expression and (\ref{5.03}), we obtain (ii) (see the end of the proof of (i) of Lemma
\ref{5.7}).
\medskip
\end{proof*}
 
Now, we are able to prove Theorem \ref{borne+}.
\medskip \\
\begin{vproof}{of Theorem \ref{borne+}.}
As for the proof of Theorem \ref{commut1}, it is sufficient to show that \linebreak $|\:\|xu(t)\|_{L^2} - \| xT(t)u_+\|_{L^2}|$
remains bounded as $t\tends\infty.$ Lemmas \ref{5.7} and \ref{5.8} achieve the proof and give the desired estimates.
\medskip
\end{vproof}
 
Now, we establish 2 lemmas which will be used to prove Theorem \ref{borne-}. The proof is very similar to the Theorem 
\ref{borne+} one.

\begin{lem}
\label{5.7'}
Let $u$ satisfying $(\ref{hypo})$ and let $h$ be defined by $(\ref{h}).$ Assume that $\lambda>0.$ Then the following
holds. \medskip \\
$
\begin{array}{rl}
\text{\rm (i)}   & If\;\; \alpha<\dfrac{4}{N}\;\; then\; \vliminf_{t\to\infty}(\|xu(t)\|_{L^2}-\|xT(t)u_+\|_{L^2})
                                 \ge\dfrac{h'(0)+4(xu_+,i\nabla u_+)}{4\|\nabla u_+\|_{L^2}}. \medskip \\
\text{\rm (ii)}  & If\;\; \alpha>\dfrac{4}{N}\;\; then\; \vlimsup_{t\to\infty}(\|xu(t)\|_{L^2}-\|xT(t)u_+\|_{L^2})
                                 \le\dfrac{h'(0)+4(xu_+,i\nabla u_+)}{4\|\nabla u_+\|_{L^2}}. \medskip \\
\text{\rm (iii)} & If\;\; \alpha=\dfrac{4}{N}\;\; then\; \vlim_{t\to\infty}(\|xu(t)\|_{L^2}-\|xT(t)u_+\|_{L^2})
                                =\dfrac{h'(0)+4(xu_+,i\nabla u_+)}{4\| \nabla u_+\|_{L^2}}.
\end{array}
$
\end{lem}
\bigskip

\begin{proof*}
We proceed as for the proof of Lemma \ref{5.7}, using (\ref{5.51'}) instead of (\ref{5.51}).
\medskip
\end{proof*}

\begin{lem}
\label{5.8'}
Let $u$ satisfying $(\ref{hypo}).$ Assume that $\lambda>0.$ Then the following holds. \medskip \\
$
\begin{array}{rl}
 ({\rm i}) & If\;\; \alpha <\dfrac{4}{N}\;\; then\; \vlimsup_{t\to\infty}(\|xu(t)\|_{L^2}-\|xT(t)u_+\|_{L^2}) \le
             \dfrac{C}{\|\nabla u_\pm\|_{L^2}}. \medskip \\
({\rm ii}) & If\;\; \alpha >\dfrac{4}{N}\;\; then\; \vliminf_{t\to\infty}(\|xu(t)\|_{L^2}-\|xT(t)u_+ \|_{L^2})
             \ge -\dfrac{C}{\|\nabla u_\pm\|_{L^2}}.
\end{array}
$
\medskip \\
Furthermore, $C=C(\vsup_{t\in\R}\|T(-t)u(t)\|_X,N,\alpha,\lambda).$
\medskip

\end{lem}
\begin{proof*} We proceed in 2 steps.\\
{\bf Step 1.} $\forall t\ge 1$, $-\dsp\vint_{1}^{t}\dsp\vint_{1}^{s}\|\nabla u(\sigma)\|_{L^2}^2d\sigma ds
\ge -C-Ct^{-\frac{N\alpha -4}{2}}-Ct+\|\nabla u_+\|_{L^2}^2t-\dfrac{1}{2}\|\nabla u_+\|_{L^2}^2t^2,$ where
$C=C(\vsup_{t\in\R}\|T(-t)u(t)\|_X,N,\alpha,\lambda).$ \\
By Proposition \ref{propdecayX}, we have $\|u(\sigma)\|_{L^{\alpha+2}}^{\alpha+2}\le
C(\vsup_{t\in\R}\|T(-t)u(t)\|_X,N,\alpha)\sigma^{-\frac{N\alpha}{2}},$ for all $\sigma >0.$ With the  conservation of energy and
the equality $\| \nabla u_+\|_{L^2}^2=2E(\vphi)$, we deduce that for all $\sigma >0,$
$$
\| \nabla u_+\|_{L^2}^{2} -\|\nabla u(\sigma)\|_{L^2}^2 \ge
-C(\vsup_{t\in\R}\|T(-t)u(t)\|_X,N,\alpha,\lambda)\sigma^{-\frac{N\alpha}{2}}
$$
(since $\lambda>0).$ Integrating this expression over $[1,t]\times [1,s]$, we get the desired result.
\\
{\bf Step 2.} Conclusion.\\
First, note that we have for all $t>1,$
$$
\begin{array}{rl}
    & \dsp\int_0^t\int_0^s\|\nabla u(\sigma)\|_{L^2}^2d\sigma ds \\
  = & \dsp\int_0^t\int_0^1\|\nabla u(\sigma)\|_{L^2}^2d\sigma ds + \int_0^1\int_1^s\|\nabla u(\sigma)\|_{L^2}^2d\sigma ds
                                                                 + \int_1^t\int_1^s\|\nabla u(\sigma)\|_{L^2}^2d\sigma ds \\
\le & \dsp\int_0^t\int_0^1\|\nabla u(\sigma)\|_{L^2}^2d\sigma ds + \int_1^t\int_1^s\|\nabla u(\sigma)\|_{L^2}^2d\sigma ds.
\end{array}
$$
And so,
\begin{gather*}
\label{5.85'}
\vint_0^t\vint_0^s\|\nabla u(\sigma)\|_{L^2}^2d\sigma ds \le C(\|\vphi\|_{H^1},N,\alpha,\lambda)t
             + \vint_1^t\vint_1^s\|\nabla u(\sigma)\|_{L^2}^2d\sigma ds,
\end{gather*}
for all $t>1.$ And we proceed as for the proof of Lemma \ref{5.8}.
\medskip
\end{proof*}

Now, we are able to prove Theorem \ref{borne-}.
\medskip

\begin{vproof}{of Theorem \ref{borne-}.}
As for the proof of Theorem \ref{commut1}, it is sufficient to show that \linebreak $|\:\|xu(t)\|_{L^2} - \| xT(t)u_+\|_{L^2}|$
remains bounded as $t\tends\infty.$ Lemmas \ref{5.7'} and \ref{5.8'} achieve the proof and give the desired estimates.
\medskip
\end{vproof}

\begin{vproof}{of Proposition \ref{autreborne}.}
By Cauchy-Schwarz' inequality, we have
\begin{gather}
\label{autreborne01}
|\:\|xu(t)\|_{L^2}-2t\|\nabla u(t)\|_{L^2}|\le\|xT(-t)u(t)\|_{L^2},
\end{gather}
for all $t\in\R$ and for every $\lambda\not=0.$ \\
We have also the following estimate.
\begin{gather}
\label{autreborne02}
\|xT(t)u_+\|_{L^2}=2t\|\nabla u_+\|_{L^2}-\dfrac{(xu_+,i\nabla u_+)}{\|\nabla u_+\|_{L^2}}+\dfrac{\|xu_+\|_{L^2}^2}{4\|\nabla
                   u_+\|_{L^2}t}+o\left(\dfrac{1}{t}\right),
\end{gather}
for all $t>0$ large enough and for every $\lambda\not=0.$ \\
We first establish \ref{autreborne1} in the case $t>0$ that we note in this proof $\ref{autreborne1}_+.$ \ref{autreborne1} in the
case $t<0$ is obviously noted $\ref{autreborne1}_-.$ By (\ref{autreborne01}), (\ref{autreborne02}) and (\ref{5.51}), we have
$$
\begin{array}{rl}
    & \|xu(t)\|_{L^2}-\|xT(t)u_+\|_{L^2} \medskip \\
\le & 2t\|\nabla u(t)\|_{L^2}+\|xT(-t)u(t)\|_{L^2}-2t\|\nabla u_+\|_{L^2}
      +\dfrac{(xu_+,i\nabla u_+)}{\|\nabla u_+\|_{L^2}}
           -\dfrac{\|xu_+\|_{L^2}^2}{4\|\nabla u_+\|_{L^2}t}+o\left(\dfrac{1}{t}\right) \medskip \\
\le & \|xT(-t)u(t)\|_{L^2}+\dfrac{(xu_+,i\nabla u_+)}{\|\nabla u_+\|_{L^2}}-\dfrac{\|xu_+\|_{L^2}^2}{4\|\nabla
           u_+\|_{L^2}t}+o\left(\dfrac{1}{t}\right)
\end{array}
$$
for all $t>0$ large enough. \\
Hence $\ref{autreborne1}_+$ by taking $\vlimsup_{t\to\infty}$ in the above expression. $\ref{autreborne1}_-$ follows by applying
$\ref{autreborne1}_+$ to $v$ solution of (\ref{nls}) such that $v(0)=\ovl{\vphi}$. Indeed, by uniqueness $v(t)=\ovl{u(-t)}$
and so $v_+=\ovl{u_-}$. Then, using $\ref{autreborne1}_+$ and the identity $T(t)\ovl{\psi}=\ovl{T(-t)\psi}$ which holds for all
$t\in\R$ and $\psi\in L^2,$ we obtain $\ref{autreborne1}_-.$ \\
Using (\ref{5.51'}) instead of (\ref{5.51}), we obtain $\ref{autreborne2}$ by the same way.
\end{vproof}

\section{Proof of Theorem \ref{commut2}}
\label{proof2}

\begin{vproof}{of Theorem \ref{commut2}.}
We proceed in 2 steps. \\
{\bf Step 1.} We have \ref{commut21} and \ref{commut23}. \\
It is well-known that there exist $v_\pm\in X$ such that $T(-t)u(t)\xrightarrow[t\to\pm\infty]{X}v_\pm$ 
(Cazenave and Weissler \cite{MR93d:35150}). The result follows from Remark \ref{rmkcommut3}. \\
{\bf Step 2.} We have \ref{commut22}. \\
We set $v(t,x)=(x+2it\nabla)u(t,x),$ $w(t,x)=e^{-i\frac{|x|^2}{4t}}u(t,x)$ and $f(u)=\lambda|u|^\alpha u.$ Since $T(t)$ is an
isometry on $H^1(\R^N),$ we only have to show that $xT(-t)u(t)\xrightarrow[t\to\pm\infty]{L^2(\R^N)}xu_\pm.$ We have
$X\hookrightarrow L^{\alpha+2}(\R^N),$ thus $u(t)-T(t)u_\pm\xrightarrow[t\to\pm\infty]{L^{\alpha+2}(\R^N)}0$ and so
$
\vlim_{t\to\pm\infty}\|u(t)\|_{L^{\alpha+2}}=\vlim_{t\to\pm\infty}\|T(t)u_\pm\|_{L^{\alpha+2}}=0.
$
Therefore, since $\alpha>\dfrac{4}{N},$
\begin{gather}
\label{commut221}
u\in L^q(\R;W^{1,r}(\R^N)), \\
\label{commut222}
v\in L^q(\R;L^r(\R^N)),
\end{gather}
for every admissible pair $(q,r)$ (see for example Cazenave \cite{caz1}, Theorem 7.5.3 for (\ref{commut221}); following the
proof of this theorem with $v$ instead of $u,$ yields (\ref{commut222})).\\
From (\ref{commut222}) we have in particular, $v\in L^\infty(\R;L^2(\R^N))$ and so by Proposition \ref{propdecayX},
\ref{propdecayX1},
\begin{gather}
\label{commut223}
\|u(t)\|_{L^r}\le C|t|^{-\frac{2}{q}},
\end{gather}
for every admissible pair $(q,r)$ and for all $t\not=0.$ \\
We have the following integral equation. For all $t\in\R,$
$$
u(t)=T(t)u_\pm-i\vint_t^{\pm\infty}T(t-s)f(u(s))ds,
$$
from which we deduce,
\begin{gather}
\label{commut224}
\forall t\in\R,\;
T(t)(xT(-t)u(t)-xu_\pm)=-i\vint_t^{\pm\infty}T(t-s)(x+2is\nabla)f(u(s))ds.
\end{gather}
We also have $(x+2it\nabla)u(t,x)=2ite^{i\frac{|x|^2}{4t}}\nabla w(t,x).$ Moreover,
$e^{-i\frac{|x|^2}{4t}}f(u(t,x))=f(w(t,x)).$ Thus,
$
(x+2it\nabla)f(u(t,x))=2ite^{i\frac{|x|^2}{4t}}\nabla f(w(t,x)),
$
and so $|(x+2it\nabla)f(u(t,x))|=2|t||\nabla f(w(t,x))|.$ Finally,
$
\|(x+2it\nabla)f(u)\|_{L^{r'}}=2|t|\:\|\nabla f(w)\|_{L^{r'}}\le C|t|\:\||w|^\alpha\nabla w\|_{L^{r'}}.
$
From this inequality, by using the Hölder's inequality twice, we deduce that (note that $|w|=|u|$)
\begin{gather}
\label{estimv}
\|(x+2it\nabla)f(u)\|_{L^{q'}(I,L^{r'})}\le C
\|u\|_{L^{\frac{q\alpha}{q-2}}(I,L^{\frac{r\alpha}{r-2}})}^\alpha\|(x+2it\nabla)u\|_{L^q(I,L^r)},
\end{gather}
for every admissible pair $(q,r)$ and for any interval $I\subseteq\R.$
From (\ref{commut224}), from the Strichartz' estimates, from (\ref{estimv}) and from (\ref{commut222}), we have
\begin{align*}
\|xT(-t)u(t)-xu_+\|_{L^2} & = \|T(t)(xT(-t)u(t)-xu_+)\|_{L^2} \\
                          & \le C \|(x+2is\nabla)f(u)\|_{L^{q'}((t,\infty);L^{r'})} \\
                          & \le C \|u\|_{L^{\frac{q\alpha}{q-2}}((t,\infty);L^{\frac{r\alpha}{r-2}})}^\alpha\|v\|_{L^q(\R;L^r)} \\
                          & \le C \|u\|_{L^{\frac{q\alpha}{q-2}}((t,\infty);L^{\frac{r\alpha}{r-2}})}^\alpha,
\end{align*}
for every admissible pair $(q,r)$ and for all $t>0.$ \\
We use this inequality with the  admissible pair $(q,r)$ such that $r=\alpha+2$ and we apply (\ref{commut223}). Then,
$$
\|xT(-t)u(t)-xu_+\|_{L^2}\le C \left(\int_t^\infty s^{-\frac{2\alpha}{q-2}}ds\right)^{\frac{q-2}{q}}
                         \le C t^{-\frac{2\alpha-(q-2)}{q}}\xrightarrow{t\to\infty}0,
$$
since $\alpha>\dfrac{4}{N}\Longrightarrow 2\alpha>q-2.$ Hence the result. The case $t<0$ follows with the same method.
\bigskip
\end{vproof}

\baselineskip .5cm

\noindent
{\large\bf Acknowledgments} \\
I would like to express my great gratitude to Professor Thierry Cazenave, my thesis adviser, for having suggested me this work and
for encouragement and helpful advices.

\baselineskip .0cm

\bibliographystyle{abbrv}
\bibliography{BiblioPaper1}

\begin{thebibliography}{10}

\bibitem{MR86a:35121}
J.~E. Barab.
\newblock Nonexistence of asymptotically free solutions for a nonlinear
  {S}chr\"odinger equation.
\newblock {\em J. Math. Phys.}, 25(11):3270--3273, 1984.

\bibitem{caz1}
T.~Cazenave.
\newblock {\em An introduction to nonlinear {S}chr\"odinger equations},
  volume~26 of {\em Textos de Métodos Matem{\'a}ticos}.
\newblock Instituto de Matem{\'a}tica, Universidade Federal do Rio de Janeiro,
  Rio de Janeiro, third edition, 1996.

\bibitem{MR93d:35150}
T.~Cazenave and F.~B. Weissler.
\newblock Rapidly decaying solutions of the nonlinear {S}chr\"odinger equation.
\newblock {\em Comm. Math. Phys.}, 147(1):75--100, 1992.

\bibitem{MR82c:35059}
J.~Ginibre and G.~Velo.
\newblock On a class of nonlinear {S}chr\"odinger equations. {I}{I}{I}.
  {S}pecial theories in dimensions $1$, $2$\ and $3$.
\newblock {\em Ann. Inst. H. Poincar\'e Sect. A (N.S.)}, 28(3):287--316, 1978.

\bibitem{MR82c:35057}
J.~Ginibre and G.~Velo.
\newblock On a class of nonlinear {S}chr\"odinger equations. {I}. {T}he
  {C}auchy problem, general case.
\newblock {\em J. Funct. Anal.}, 32(1):1--32, 1979.

\bibitem{MR82c:35058}
J.~Ginibre and G.~Velo.
\newblock On a class of nonlinear {S}chr\"odinger equations. {I}{I}.
  {S}cattering theory, general case.
\newblock {\em J. Funct. Anal.}, 32(1):33--71, 1979.

\bibitem{MR87b:35150}
J.~Ginibre and G.~Velo.
\newblock The global {C}auchy problem for the nonlinear {S}chr\"odinger
  equation revisited.
\newblock {\em Ann. Inst. H. Poincar\'e Anal. Non Lin\'eaire}, 2(4):309--327,
  1985.

\bibitem{MR87i:35171}
J.~Ginibre and G.~Velo.
\newblock Scattering theory in the energy space for a class of nonlinear
  {S}chr\"odinger equations.
\newblock {\em J. Math. Pures Appl. (9)}, 64(4):363--401, 1985.

\bibitem{MR88f:35133}
T.~Kato.
\newblock On nonlinear {S}chr\"odinger equations.
\newblock {\em Ann. Inst. H. Poincar\'e Phys. Th\'eor.}, 46(1):113--129, 1987.

\bibitem{MR1646048}
M.~Keel and T.~Tao.
\newblock Endpoint {S}trichartz estimates.
\newblock {\em Amer. J. Math.}, 120(5):955--980, 1998.

\bibitem{MR1726753}
K.~Nakanishi.
\newblock Energy scattering for nonlinear {K}lein-{G}ordon and {S}chr\"odinger
  equations in spatial dimensions $1$ and $2$.
\newblock {\em J. Funct. Anal.}, 169(1):201--225, 1999.

\bibitem{MR1829982}
K.~Nakanishi.
\newblock Remarks on the energy scattering for nonlinear {K}lein-{G}ordon and
  {S}chr\"odinger equations.
\newblock {\em Tohoku Math. J. (2)}, 53(2):285--303, 2001.

\bibitem{297.35062}
W.~A. Strauss.
\newblock Nonlinear scattering theory.
\newblock In J.~Lavita and J.-P. Planchard, editors, {\em Scattering Theory in
  Mathematical Physics}, pages 53--78, Reidel, 1974.

\bibitem{MR83b:47074a}
W.~A. Strauss.
\newblock Nonlinear scattering theory at low energy.
\newblock {\em J. Funct. Anal.}, 41(1):110--133, 1981.

\bibitem{MR87g:35221}
Y.~Tsutsumi.
\newblock Scattering problem for nonlinear {S}chr\"odinger equations.
\newblock {\em Ann. Inst. H. Poincar\'e Phys. Th\'eor.}, 43(3):321--347, 1985.

\bibitem{MR85k:35216}
Y.~Tsutsumi and K.~Yajima.
\newblock The asymptotic behavior of nonlinear {S}chr\"odinger equations.
\newblock {\em Bull. Amer. Math. Soc. (N.S.)}, 11(1):186--188, 1984.

\end{thebibliography}

\end{document}